%% file: OPMWGraphs.tex
\documentclass[oneside,english]{amsart}
\usepackage[latin1]{inputenc}
\usepackage{amssymb}

\makeatletter
 \theoremstyle{plain}    
 \newtheorem{thm}{Theorem}[section]
 \numberwithin{equation}{section} 
 \numberwithin{figure}{section} 
 \theoremstyle{plain}
 \theoremstyle{plain}    
 \newtheorem{prop}[thm]{Proposition} 
 \theoremstyle{definition}
 \newtheorem{defn}[thm]{Definition}
 \theoremstyle{plain}    
 \newtheorem{lem}[thm]{Lemma} 
 \theoremstyle{plain}    
 \newtheorem{cor}[thm]{Corollary} 
 \theoremstyle{definition}
  \newtheorem{example}[thm]{Example}

\usepackage[dvips]{graphicx}
\usepackage[dvips]{color}

\usepackage{amsmath}
\usepackage{amsfonts}%
\usepackage{amssymb}

\usepackage{babel}
\makeatother
\begin{document}
\address[Marius Ionescu]{Department of Mathematics, The University of Iowa, Iowa City, Iowa, USA, email: mionescu@math.uiowa.edu}
\address[Yasuo Watatani]{Department of Mathematical Sciences,  Kyushu University, Hakozaki,  Fukuoka, 812-8581,  Japan, \\ email: watatani@math.kyushu-u.ac.jp}

\title{$C^{\ast}$-Algebras associated with Mauldin-Williams Graphs}

\author{Marius Ionescu}

\author{and Yasuo Watatani }

\begin{abstract}
A Mauldin-Williams graph $\mathcal{M}$ is a generalization of an
iterated function system by a directed graph. Its invariant set $K$
plays the role of the self-similar set. We associate a $C^{*}$-algebra
$\mathcal{O}_{\mathcal{M}}(K)$ with a Mauldin-Williams graph $\mathcal{M}$
and the invariant set $K$ laying emphasis on the singular points.
We assume that the underlying graph $G$ has no sinks and no sources.
If $\mathcal{M}$ satisfies the open set condition in $K$ and $G$
is irreducible and is not a cyclic permutation, then the associated
$C^{*}$-algebra $\mathcal{O}_{\mathcal{M}}(K)$ is simple and purely
infinite. We calculate the K-groups for some examples including the
inflation rule of the Penrose tilings. 
\end{abstract}
\maketitle

\section{Introduction}

Self-similar sets are often constructed as the invariant set of iterated
function systems \cite{Barnsley}, \cite{Hutchinson}, \cite{kigami}.
Many other examples, such as the inflation rule of the Penrose tilings,
show that a graph directed generalization of iterated function systems
is also interesting and has been developed as Mauldin-Williams graphs
\cite{BarGrigNek}, \cite{Edgar-Mauldin}, \cite{Mauldin-Williams}.
Since we can regard graph directed iterated function systems as dynamical
systems, we expect that there exist fruitful connections between Mauldin-Williams
graphs and $C^{*}$-algebras. 

In \cite{John} the first-named author defined a $C^{*}$-correspondence
$\mathcal{X}$ for a Mauldin-Williams graph $\mathcal{M}$ and showed
that the Cuntz-Pimsner algebra $\mathcal{O}_{\mathcal{X}}$ (see \cite{Fowler-Muhly-Raeburn},
\cite{Pimnser}) is isomorphic to the Cuntz-Krieger algebra $\mathcal{O}_{G}$
(\cite{Cuntz-Krieger}) associated with the underlying graph $G$.
In particular, if the Mauldin-Williams graph has one vertex and $N$
edges, the Cuntz-Pimsner algebra $\mathcal{O}_{\mathcal{X}}$ is isomorphic
to the Cuntz algebra $\mathcal{O}_{N}$ (\cite{Cuntz}), which recovers
a result in \cite{Pinzari-Watatami-Y}. The construction is useful
because it gives many examples of different non-self-adjoint algebras
which are not completely isometrically isomorphic, but have the same
$C^{*}$-envelope (\cite{Muhly-Baruc-tensoalgebras}) as shown by
the first-named author in \cite{Io}. 

On the other hand, Kajiwara and the second-named author introduced
$C^{*}$-algebras associated with rational functions in \cite{KW}
including the singular points (i.e. branched points). They developed
the idea to associate $C^{*}$-algebras with self-similar sets considering
the singular points in \cite{KWp03}. In this paper we associate another
$C^{*}$-algebra $\mathcal{O}_{\mathcal{M}}(K)$ with a Mauldin-Williams
graph $\mathcal{M}$ and its invariant set $K$ putting emphasis on
the singular points as above. We show that the associated $C^{*}$-algebra
$\mathcal{O}_{\mathcal{M}}(K)$ is not isomorphic to the Cuntz-Krieger
algebra $\mathcal{O}_{G}$ for the underlying graph $G$ in general.
This comes from the fact that the singular points cause the failure
of the injectivity of the coding by the Markov shift for $G$. We
assume that the underlying graph G has no sinks and no sources. We
show that the associated $C^{*}$-algebra $\mathcal{O}_{\mathcal{M}}(K)$
is simple and purely infinite if $\mathcal{M}$ satisfies the open
set condition in $K$ and if $G$ is irreducible and is not a cyclic
permutation. We calculate the $K$-groups for some examples including
the inflation rule of the Penrose tilings. The $C^{*}$-algebras associated
with tilings were first studied by Connes in \cite{Co} and also discussed
by Mingo \cite{Mi} and Anderson-Putnam \cite{AP}. But we do not
know the exact relatition between their constructions and ours.

Our construction has a common flavor with several topological generalizations
of graph $C^{*}$-algebras \cite{Kumjian-Pask-Raeburn-directedgraphs}
by Brenken \cite{Br}, Katsura \cite{Ka}, Muhly and Solel \cite{Muhly-Baruc-Moritaequiv}
and Muhly and Tomforde \cite{MT}. Since graph directed iterated function
systems are sometimes obtained as continuous cross sections of expanding
maps, $C^{*}$-algebras associated with interval maps introduced by
Deaconu and Shultz in \cite{DS} are closely related with our construction.
In a different point of view, Bratteli and Jorgensen studied a relation
between iterated function systems and representation of Cuntz algebras
in \cite{BJ}. 

The authors would like to thank P. Muhly for suggesting their collaboration
and giving them the chance to discuss this work at a conference at
the University of Iowa.

\section{Mauldin-Williams graphs and the associated $C^{\ast}$-correspondence}

By a \emph{Mauldin-Williams} \emph{graph} we mean a system $\mathcal{M}=(G,\{ T_{v},\rho_{v}\}_{v\in E^{0}},\{\phi_{e}\}_{e\in E^{1}})$,
where $G=(E^{0},E^{1},r,s)$ is a graph with a finite set of vertices
$E^{0}$, a finite set of edges $E^{1}$, a \emph{range} map $r$
and a \emph{source} map $s$, and where $\{ T_{v},\rho_{v}\}_{v\in E^{0}}$
and $\{\phi_{e}\}_{e\in E^{1}}$ are families such that:

\begin{enumerate}
\item Each $T_{v}$ is a compact metric space with a prescribed metric $\rho_{v}$,
$v\in E^{0}$.
\item For $e\in E$, $\phi_{e}$ is a continuous map from $T_{r(e)}$ to
$T_{s(e)}$ such that\[
c^{\prime}\rho_{r(e)}(x,y)\le\rho_{s(e)}(\phi_{e}(x),\phi_{e}(y))\leq c\rho_{r(e)}(x,y)\]
 for some constants $c^{\prime},c$ satisfying $0<c^{\prime}<c<1$
(independent of $e$) and all $x,y\in T_{r(e)}$. 
\end{enumerate}
We shall assume, too, that the source map $s$ and the range map $r$
are surjective. Thus, we assume that there are no sinks and no sources
in the graph $G$.
\newcommand{\twolines}[2]{{\genfrac{}{}{0pt}{}{#1}{#2}}}

In the particular case when we have only one vertex and $N$ edges,
we obtain a so called \emph{iterated function system} $(K,\{\phi_{i}\}_{i=1,\dots,N})$.

An \emph{invariant list} associated with a Mauldin-Williams graph
$\mathcal{M}=(G,\{ T_{v},\rho_{v}\}_{v\in E^{0}},\{\phi_{e}\}_{e\in E^{1}})$
is a family $(K_{v})_{v\in E^{0}}$ of compact sets such that $K_{v}\subset T_{v}$
for all $v\in E^{0}$, and such that\begin{align*}
K_{v} & =\bigcup_{\twolines{e\in E^{1}}{s(e)=v}}\phi_{e}(K_{r(e)}).\end{align*}
Since each $\phi_{e}$ is a contraction, $\mathcal{M}$ has a unique
invariant list (see \cite[Theorem 1]{Mauldin-Williams}). We set $T:=\bigcup_{v\in E^{0}}T_{v}$
and $K:=\bigcup_{v\in E^{0}}K_{v}$ and we call $K$ the \emph{invariant
set} of the Mauldin-Williams graph.

In the particular case when we have one vertex $v$ and $n$ edges,
i.e. in the setting of an \emph{iterated function system,} the invariant
set is the unique compact subset $K:=K_{v}$ of $T=T_{v}$ such that

\newcommand{\MWG}{\mathcal{M}=(G,\{ K_{v},\rho_{v}\}_{v\in E^{0}},\{\phi_{e}\}_{e\in E^{1}})}
\[
K=\phi_{1}(K)\cup\cdots\cup\phi_{n}(K).\]
In this paper we forget about the ambient space $T$. That is, we
consider that the Mauldin-Williams graph is $\MWG$, with $K=\bigcup_{v\in E^{0}}K_{v}$
being the invariant set.

We say that a graph $G=(E^{0},E^{1},r,s)$ is \emph{irreducible} (or
\emph{totally connected}) if for every $v_{1},v_{2}\in E^{0}$ there
is a finite path $w$ such that $s(w)=v_{1}$ and $r(w)=v_{2}$. We
will assume in this paper that the graph $G=(E^{0},E^{1},r,s)$ is
irreducible and not a cyclic permutation. That is there exists a vertex
$v^{\ast}\in E^{0}$ and there exist two edges $e_{1}\ne e_{2}$ such
that $s(e_{1})=s(e_{2})=v^{\ast}$.

For a natural number $m$, we define $E^{m}:=\{ w=(w_{1},\dots,w_{m})\::\: w_{i}\in E^{1}\;\mbox{and }r(w_{i})=s(w_{i+1})\mbox{ for all }i=1,\dots,m-1\}$.
An element $w\in E^{m}$ is called a path of length $m$. We extend
$s$ and $r$ to $E^{m}$ by $s(w)=s(w_{1})$ and $r(w)=r(w_{m})$
for all $w\in E^{m}$. We set $E^{*}=\bigcup_{m\ge1}E^{m}$ and denote
the length of a path $w$ by $l(w)$. We define also $E^{m}(v):=\{ w\in E^{m}\;:\; s(w)=v\}$
to be the set of paths of length $m$ starting at the vertex $v$
and $E^{\ast}(v)=\bigcup_{m\ge1}E^{m}(v)$ to be the set of finite
paths starting at $v$.

The infinite path space is $E^{\infty}=\{ w=(w_{n})_{n\ge1}\::\: r(w_{n})=s(w_{n+1})\mbox{ for all }n\ge1\}$
and the space of infinite paths starting at a vertex $v$ is $E^{\infty}(v)=\{ w\in E^{\infty}\;:\; s(w)=s(w_{1})=v\}$.
On $E^{\infty}(v)$ we define the metric $\delta_{v}(\alpha,\beta)=c^{l(\alpha\wedge\beta)}$
if $\alpha\ne\beta$ and $0$ otherwise, where $\alpha\wedge\beta$
is the longest common prefix of $\alpha$ and $\beta$ (see \cite[Page 116]{Edgar}).
Then $E^{\infty}(v)$ is a compact metric space, and, since $E^{\infty}$
equals the disjoint union of the spaces $E^{\infty}(v)$, $E^{\infty}$
is a compact metric space in a natural way. For $w\in E^{m}$, let
$\phi_{w}=\phi_{w_{1}}\circ\cdots\phi_{w_{m}}$ and $K_{w}=\phi_{w}(K_{r(w)})$.
Then for any infinite path $\alpha=(\alpha_{n})_{n\in\mathbb{N}}$,
$\bigcap_{n\ge1}K_{(\alpha_{1},\dots,\alpha_{n})}$ contains only
one point. Therefore we can define a map $\pi:E^{\infty}\rightarrow K$
by $\{\pi(\alpha)\}=\bigcap_{n\ge1}K_{(\alpha_{1},\dots,\alpha_{n})}$.
Since $\pi(E^{\infty})$ is also an invariant set, we have that $\pi(E^{\infty})=K$.
Thus $\pi$ is a continuous, onto map. Moreover, for any $y\in K_{v_{0}}$
and any neighborhood $U\subset K_{v_{0}}$ of $y$, there exist $n\in\mathbb{N}$
and $w\in E^{n}(v_{0})$ such that\[
y\in\phi_{w}(K_{r(w)})\subset U.\]

We say that a Mauldin-Williams graph $\mathcal{M}$ satisfies the
\emph{open set condition in $K$} if there exists a family of non-empty
sets $(V_{v})_{v\in E^{0}}$, such that $V_{v}\subset K_{v}$ for
all $v\in E^{0}$, and such that\[
\bigcup_{\twolines{e\in E^{1}}{s(e)=v}}\phi_{e}(V_{r(e)})\subset V_{v}\;\mbox{for all }v\in E^{0},\]
\[
\mbox{and }\phi_{e}(V_{r(e)})\bigcap\phi_{f}(V_{r(f)})=\emptyset\;\textrm{\mbox{if }}e\ne f.\]
Then $V:=\bigcup_{v\in E^{0}}V_{v}$ is an open dense subset of $K$.
Moreover, for $n\in\mathbb{N}$ and $w,v\in E^{n}$, if $w\ne v$
and $r(w)=r(v)$ then $\phi_{w}(V_{r(w)})\bigcap\phi_{v}(V_{r(v)})=\emptyset$.

For $e\in E^{1}$, we define the \emph{cograph of} $\phi_{e}$ to
be the set 
\newcommand{\cograph}{\operatorname{cograph}}
 \[
\cograph(\phi_{e})=\{(x,y)\in K_{s(e)}\times K_{r(e)}\;:\; x=\phi_{e}(y)\}\subset K\times K.\]
We shall consider the union\[
\mathcal{G}=\mathcal{G}(\{\phi_{e}\;:\: e\in E^{1}\}):=\bigcup_{e\in E^{1}}\cograph(\phi_{e}).\]
Consider the $C^{\ast}$-algebra $A=C(K)$ and let $X=C(\mathcal{G})$.
Then $X$ is a $C^{\ast}$-correspondence over $A$ with the structure
defined by the formulae:\begin{eqnarray*}
(a\cdot\xi\cdot b)(x,y) & = & a(x)\xi(x,y)b(y),\\
\langle\xi,\eta\rangle_{A}(y) & = & \sum_{\twolines{e\in E^{1}}{y\in K_{r(e)}}}\overline{\xi(\phi_{e}(y),y)}\eta(\phi_{e}(y),y),\end{eqnarray*}
for all $a,b\in A$, $\xi,\eta\in X$, $(x,y)\in\mathcal{G}$ and
$y\in K$. It is clear that the $A$-valued inner product is well
defined. The left multiplication is given by the $\ast$-homomorphism
$\Phi:A\rightarrow\mathcal{L}(X)$ such that $(\Phi(a)\xi)(x,y)=a(x)\xi(x,y)$
for $a\in A$ and $\xi\in X$. Put $\Vert\xi\Vert_{2}=\Vert\langle\xi,\xi\rangle_{A}\Vert_{\infty}^{1/2}$.

For any natural number $n$, we define $\mathcal{G}_{n}=\mathcal{G}(\{\phi_{w}\::\: w\in E^{n}\})$
and a $C^{\ast}$-correspondence $X_{n}=C(\mathcal{G}_{n})$ similarly.
We also define a path space $\mathcal{P}_{n}$ of length $n$ by\begin{eqnarray*}
\mathcal{P}_{n} & = & \{(\phi_{w_{1},\cdots,w_{n}}(y),\phi_{w_{2},\cdots,w_{n}}(y),\dots,\phi_{w_{n}}(y),y)\in K^{n+1}:\\
 &  & w=(w_{1},\dots,w_{n})\in E^{n},y\in K_{r(w)}\}.\end{eqnarray*}
Then $Y_{n}:=C(\mathcal{P}_{n})$ is a $C^{\ast}$-correspondence
over $A$ with an $A$-valued inner product defined by\[
\langle\xi,\eta\rangle_{A}(y)=\sum_{\twolines{w\in E^{n}}{y\in K_{r(w)}}}\overline{\xi(\phi_{w_{1},\dots,w_{n}}(y),\dots,\phi_{w_{n}}(y),y)}\eta(\phi_{w_{1},\dots,w_{n}}(y),\dots,\phi_{w_{n}}(y),y),\]
for all $\xi,\eta\in Y_{n}$ and $y\in K$.

\begin{prop}
Let $\MWG$ be a Mauldin-Williams graph. Let $K$ be the invariant
set. Then $X=C(\mathcal{G})$ is a full $C^{\ast}$-correspondence
over $A=C(K)$ without completion. The left action $\Phi:A\rightarrow\mathcal{L}(X)$
is unital and faithful. Similar statements hold for $Y_{n}=C(\mathcal{P}_{n})$.
\end{prop}
\begin{proof}
For any $\xi\in X$ we have\[
\Vert\xi\Vert_{\infty}\le\Vert\xi\Vert_{2}=\left(\sup_{y\in K}\sum_{\twolines{e\in E^{1}}{y\in K_{r(e)}}}|\xi(\phi_{e}(y),y)|^{2}\right)^{1/2}\le\sqrt{N}\Vert\xi\Vert_{\infty},\]
where $N$ is the number of edges in $E^{1}$. Therefore the norms
$\Vert\;\Vert_{2}$ and $\Vert\;\Vert_{\infty}$ are equivalent. Since
$C(\mathcal{G})$ is complete with respect to $\Vert\;\Vert_{\infty}$,
it is also complete with respect to $\Vert\;\Vert_{2}$.

Let $\xi\in X$ be defined by the formula\[
\xi(x,y)=\frac{1}{\sqrt{\#(e:y\in K_{r(e)})}}\quad\mbox{for all }(x,y)\in\mathcal{G}.\]
Then $\langle\xi,\xi\rangle_{A}(y)=1$, hence $\langle X,X\rangle_{A}$
contains the identity of $A$. Therefore $X$ is full. If $a\in A$
is not zero, then there exists $x_{0}\in K$ such that $a(x_{0})\ne0$.
Since $K$ is the invariant set of the Mauldin-Williams graph, there
exists $e\in E^{1}$ and $y_{0}\in K_{r(e)}$ such that $x_{0}\in K_{s(e)}$
and $\phi_{e}(y_{0})=x_{0}$. Choose $\xi\in X$ such that $\xi(x_{0},y_{0})\ne0$.
Then $\Phi(a)\xi\ne0$, hence $\Phi$ is faithful. The statements
for $Y_{n}$ are similarly proved.
\end{proof}
\begin{defn}
Let $\MWG$ be a Mauldin-Williams graph with the invariant set $K$.
We associate a $C^{\ast}$-algebra $\mathcal{O}_{\mathcal{M}}(K)$
to $\mathcal{M}$ as the Cuntz-Pimsner algebra $\mathcal{O}_{X}$
of the $C^{\ast}$-correspondence $X=C(\mathcal{G})$ over the $C^{\ast}$-algebra
$A=C(K)$.
\end{defn}
As in \cite{KWp03}, we denote by $\mathcal{O}_{X}^{alg}$ the $\ast$-algebra
generated algebraically by $A$ and $S_{\xi}$ with $\xi\in X$. The
\emph{gauge action} is $\gamma:\mathbb{R}\rightarrow\operatorname{Aut}\mathcal{O}_{X}$
defined by $\gamma_{t}(S_{\xi})=e^{it}S_{\xi}$ for all $\xi\in X$,
and $\gamma_{t}(a)=a$ for all $a\in A$.

\begin{prop}
\label{pro:Embeddings}Let $\MWG$ be a Mauldin-Williams graph. Assume
that $K$ is the invariant set of the graph. Then there is an isomorphism
$\varphi_{n}:X^{\otimes n}\rightarrow C(\mathcal{P}_{n})$ as $C^{\ast}$-correspondences
over $A$ such that\[
(\varphi_{n}(\xi_{1}\otimes\cdots\otimes\xi_{n}))(\phi_{w_{1},\dots,w_{n}}(y),\dots,\phi_{w_{n}}(y),y)\]
\[
=\xi_{1}(\phi_{w_{1},\dots,w_{n}}(y),\phi_{w_{2},\dots,w_{n}}(y))\xi_{2}(\phi_{w_{2},\dots,w_{n}}(y),\phi_{w_{3},\dots,w_{n}}(y))\dots\xi_{n}(\phi_{w_{n}}(y),y)\]
for all $\xi_{1},\dots,\xi_{n}\in X$, $y\in K$, and $w=(w_{1},\dots,w_{n})\in E^{n}$
such that $y\in K_{r(w_{n})}$. Moreover, let $\rho_{n}:\mathcal{P}_{n}\rightarrow\mathcal{G}_{n}$
be an onto continuous map such that\[
\rho_{n}(\phi_{w_{1},\dots,w_{n}}(y),\dots,\phi_{w_{n}}(y),y)=(\phi_{w_{1},\dots,w_{n}}(y),y).\]
Then $\rho_{n}^{\ast}:C(\mathcal{G}_{n})\ni f\rightarrow f\circ\rho_{n}\in C(\mathcal{P}_{n})$
is an embedding as a Hilbert submodule preserving inner product.
\end{prop}
\begin{proof}
It is easy to see that $\varphi_{n}$ is well-defined and a bimodule
morphism. We show that $\varphi_{n}$ preserves inner product. Consider
the case when $n=2$ for simplicity of the notation. Let $\xi_{1},\xi_{2},\eta_{1},\eta_{2}\in X$.
We have\begin{eqnarray*}
\langle\xi_{1}\otimes\xi_{2},\eta_{1}\otimes\eta_{2}\rangle_{A}(y) & = & \langle\xi_{2},\langle\xi_{1},\eta_{1}\rangle_{A}\cdot\eta_{2}\rangle_{A}(y)\\
 & = & \sum_{{\genfrac{}{}{0pt}{}{e\in E}{y\in K_{r(e)}}}}\overline{\xi_{2}(\phi_{e}(y),y)}\langle\xi_{1},\eta_{1}\rangle_{A}(\phi_{e}(y))\eta_{2}(\phi_{e}(y),y)\\
 & = & \sum_{{\genfrac{}{}{0pt}{}{fe\in E^{2}}{y\in K_{r(e)}}}}\overline{\xi_{1}(\phi_{fe}(y),\phi_{e}(y))\xi_{2}(\phi_{e}(y),y)}\eta_{1}(\phi_{fe}(y),\phi_{e}(y))\eta_{2}(\phi_{e}(y),y)\\
 & = & \langle\varphi_{2}(\xi_{1}\otimes\xi_{2}),\varphi_{2}(\eta_{1}\otimes\eta_{2})\rangle_{A}(y).\end{eqnarray*}
Since $\varphi_{n}$ preserves the inner product, it is one to one.
Using the Stone-Weierstrass Theorem, one can show that $\varphi_{n}$
is also onto. The rest is clear.
\end{proof}
We let $i_{n,m}:C(\mathcal{P}_{n})\rightarrow C(\mathcal{P}_{m})$
be the natural inner-product preserving embedding, for $m\ge n$.

\begin{defn}
Consider a covering map $\pi:\mathcal{G}\rightarrow K$ defined by
$\pi(x,y)=y$ for $(x,y)\in\mathcal{G}$. Define the set\[
B(\mathcal{M}):=\{ x\in K\;:\; x=\phi_{e}(y)=\phi_{f}(y)\mbox{ for some }y\in K\mbox{ and }e\ne f\}.\]
The set $B(\mathcal{M})$ will be described by the ideal $I_{X}:=\Phi^{-1}(\mathcal{K}(X))$
of $A$. We define a branch index $e(x,y)$ at $(x,y)\in\mathcal{G}$
by\[
e(x,y):=\#\{ e\in E^{1}\;:\;\phi_{e}(y)=x\}.\]
Hence $x\in B(\mathcal{M})$ if and only if there exists some $y\in K$
with $e(x,y)\ge2$. For $x\in K$ we define\[
I(x):=\{ e\in E^{1}\;:\;\textrm{there exists }y\in K\mbox{ such that }x=\phi_{e}(y)\}.\]

\end{defn}
\begin{lem}
In the above situation, if $x\in K\setminus B(\mathcal{M})$, then
there exists an open neighborhood $U_{x}$ of $x$ satisfying the
following:
\begin{enumerate}
\item $U_{x}\bigcap B(\mathcal{M})=\emptyset$.
\item If $e\in I(x)$, then $\phi_{f}(\phi_{e}^{-1}(U_{x}))\bigcap U_{x}=\emptyset$
for $e\ne f$, such that $r(e)=r(f)$.
\item If $e\notin I(x)$, then $U_{x}\bigcap\phi_{e}(K_{r(e)})=\emptyset$.
\end{enumerate}
\end{lem}
\begin{proof}
Let $x\in K\setminus B(\mathcal{M})$. Let $v_{0}\in E^{0}$ such
that $x\in K_{v_{0}}$. Since $B(\mathcal{M})$ and $\bigcup_{e\notin I(x)}\phi_{e}(K_{r(e)})$
are closed and $x$ is not in either of them, there exists an open
neighbourhood $W_{x}\subset K_{v_{0}}$ of $x$ such that\[
W_{x}\bigcap(B(\mathcal{M})\cup\bigcup_{e\notin I(x)}\phi_{e}(K_{r(e)}))=\emptyset.\]
For each $e\in I(x)$ there exists a unique $y_{e}\in K$ with $x=\phi_{e}(y_{e})$,
since $x\notin B(\mathcal{M})$. For $f\in E^{1}$, if $r(e)=r(f)$
and $f\neq e$ then $\phi_{f}(y_{e})\neq\phi_{e}(y_{e})=x$. Therefore
there exists an open neighborhood $V_{x}^{e}$ of $x$ such that $\phi_{f}(\phi_{e}^{-1}(V_{x}^{e}))\bigcap V_{x}^{e}=\emptyset$
if $f\ne e$ and $r(f)=r(e)$. Let $U_{x}:=W_{x}\bigcap(\bigcap_{e\in I(x)}V_{x}^{e})$.
Then $U_{x}$ is an open neighborhood of $x$ and satisfies all the
requirements.
\end{proof}
\begin{prop}
Let $\MWG$ be a Mauldin-Williams graph. Assume that the system $\mathcal{M}$
satisfies the open set condition in $K$. Then\[
I_{X}=\{ a\in A=C(K)\;:\: a\mbox{ vanishes on }B(\mathcal{M})\}.\]

\end{prop}
\begin{proof}
The proof requires only minor modifications from the proof of \cite[Proposition 2.4]{KWp03}.
\end{proof}
\begin{cor}
$\# B(\mathcal{M})=\operatorname{dim}(A/I_{X})$.
\end{cor}
\;

\begin{cor}
The closed set $B(\mathcal{M})$ is empty if and only if $\Phi(A)$
is contained in $\mathcal{K}(X)$ if and only if $X$ is finitely
generated projective right $A$ module.
\end{cor}

\section{Simplicity and Pure Infinitness}

Let $\MWG$ be a Mauldin-Williams graph. Let $A=C(K)$ and $X=C(\mathcal{G})$.
For $e\in E^{1}$ define an endomorphism $\beta_{e}:A\rightarrow A$
by\[
(\beta_{e}(a))(y):=\begin{cases}
a(\phi_{e}(y)) & \textrm{\mbox{if }}y\in K_{r(e)}\\
0 & \mbox{otherwise},\end{cases}\]
for all $a\in A$ and $y\in K$. We also define a unital completely
positive map $E_{\mathcal{M}}:A\rightarrow A$ by\begin{eqnarray*}
(E_{\mathcal{M}}(a))(y): & = & \frac{1}{\#\{ e\in E^{1}:y\in K_{r(e)}\}}\sum_{\twolines{e\in E^{1}}{y\in K_{r(e)}}}a(\phi_{e}(y))\\
 & = & \frac{1}{\#\{ e\in E^{1}:y\in K_{r(e)}\}}\sum_{\twolines{e\in E^{1}}{y\in K_{r(e)}}}\beta_{e}(a)(y),\end{eqnarray*}
for $a\in A$, $y\in K$. For the function $\xi_{0}\in X$ defined
by the formula\[
\xi_{0}(x,y)=\frac{1}{\sqrt{\#\{ e\in E^{1}:y\in K_{r(e)}\}}}\]
we have\[
E_{\mathcal{M}}(a)=\langle\xi_{0},\Phi(a)\xi_{0}\rangle_{A}\;\mbox{and}\; E_{\mathcal{M}}(I)=\langle\xi_{0},\xi_{0}\rangle_{A}=I.\]
Let $D:=S_{\xi_{0}}\in\mathcal{O}_{\mathcal{M}}(K)$.

\begin{lem}
In the above situation, for $a\in A$, we have that $D^{\ast}aD=E_{\mathcal{M}}(a)$
and in particular $D^{\ast}D=I$.
\end{lem}
\begin{proof}
The same with \cite[Lemma 3.1]{KWp03}.
\end{proof}
\begin{defn}
Let $\MWG$ be a Mauldin-Williams graph. We say that an element $a\in A=C(K)$
is $(\mathcal{M},E^{n})$-invariant if\[
a(\phi_{\alpha}(y))=a(\phi_{\beta}(y))\;\mbox{for any}\; y\in K\;\mbox{and}\;\alpha,\beta\in E^{n}\;\mbox{such that}\; y\in K_{r(\alpha)}=K_{r(\beta)}.\]
If $a\in A$ is $(\mathcal{M},E^{n})$-invariant then $a$ is also
$(\mathcal{M},E^{n-1})$-invariant. Then, similar like in \cite[Definition on page 11]{KWp03},
if $a$ is $(\mathcal{M},E^{n})$-invariant, we can define\[
\beta^{k}(a)(y):=a(\phi_{w_{1}}\dots\phi_{w_{n}}(y)),\;\mbox{for any}\; w\in E^{k},\;\mbox{such that}\; y\in K_{r(w_{n})}.\]
Then, for any $\xi_{1},\dots,\xi_{n}\in X$ and $a\in A$ $(\mathcal{M},E^{n})$-invariant
, we have the relation:\[
aS_{\xi_{1}}\dots S_{\xi_{n}}=S_{\xi_{1}}\dots S_{\xi_{n}}\beta^{n}(a).\]

\end{defn}
\begin{lem}
Let $\MWG$ be a Mauldin-Williams graph such that $(E^{0},E^{1},r,s)$
is an irreducible graph. For any non-zero positive element $a\in A$
and for every $\varepsilon>0$ there exists $n\in\mathbb{N}$ and
$\xi\in X^{\otimes n}$ with $\langle\xi,\xi\rangle_{A}=1$ such that\[
\Vert a\Vert-\varepsilon\le S_{\xi}^{\ast}aS_{\xi}\le\Vert a\Vert.\]

\end{lem}
\begin{proof}
Let $x_{0}\in K$ be such that $\Vert a\Vert=a(x_{0})$. Let $v_{0}\in E^{0}$
such that $x_{0}\in K_{v_{0}}$. Then there exists an open neighborhood
$U_{0}$ of $x_{0}$ in $K_{v_{0}}$ such that for any $x\in U_{0}$
we have $\Vert a\Vert-\varepsilon\le a(x)\le\Vert a\Vert$. Let $U_{1}$
an open neighborhood of $x_{0}$ in $K_{v_{0}}$ and $K_{1}$ compact
such that $U_{1}\subset K_{1}\subset U_{0}$. Since the map $\pi:E^{\infty}\rightarrow K$
is onto and continuous, there exists some $n_{1}\in\mathbb{N}$ and
$\alpha\in E^{n_{1}}(v_{0})$ such that $\phi_{\alpha}(K_{r(\alpha)})\subset U_{1}$.
For any vertex $v\in V$, since the graph $G$ is irreducible, there
exists a path $w_{v}\in E^{\ast}$ from $r(\alpha)$ to $v$. Then
$\phi_{\alpha w_{v}}(K_{v})\subset U_{1}$. Hence, for each $v\in V$,
there is $\alpha_{v}\in E^{n(v)}(v_{0})$, for some $n(v)\ge n_{1}$,
such that $\phi_{\alpha_{v}}(K_{v})\subset U_{1}$. For each $v\in V$,
define the closed subsets $F_{1,v}$ and $F_{2,v}$ of $K\times K$
by\begin{eqnarray*}
F_{1,v} & = & \{(x,y)\in K\times K\;:\; x=\phi_{\alpha}(y),x\in K_{1},y\in K_{v},\alpha\in E^{n(v)}(v_{0})\}\\
F_{2,v} & = & \{(x,y)\in K\times K\;:\; x=\phi_{\alpha}(y),x\in U_{0}^{c},y\in K_{v},\alpha\in E^{n(v)}(v_{0})\}.\end{eqnarray*}
Since $F_{1,v}\cap F_{2,w}=\emptyset$ for all $w\in V$ and $F_{1,v}\cap F_{1,w}=\emptyset$
if $v\ne w$, there exists $g_{v}\in C(\mathcal{G}_{n(v)})$ such
that $0\le g_{v}(x,y)\le1$ and\[
g_{v}(x,y)=\begin{cases}
1 & \mbox{if }(x,y)\in F_{1,v}\\
0 & \mbox{if }(x,y)\in\cup_{w\in E^{0}}F_{2,w}\bigcup\cup_{w\ne v}F_{1,w}.\end{cases}\]
Since $\phi_{\alpha_{v}}(K_{v})\subset U_{1}$ for each $y\in K_{v}$,
there exists $x_{y}\in U_{1}$ such that $x_{y}=\phi_{\alpha_{v}}(y)\in U_{1}\subset K_{1}$.
Therefore\begin{equation}
\langle g_{v},g_{v}\rangle_{A}(y)=\sum_{\twolines{\alpha\in E^{n(v)}}{y\in K_{r(\alpha)}}}|g_{v}(\phi_{\alpha}(y),y)|^{2}\ge|g_{v}(x_{y},y)|^{2}=1\label{eq:innprodg_v}\end{equation}
for all $y\in K_{v}$. Let $n=\max\{ n(v)\;:\; v\in E^{0}\}$. We
identify $X^{\otimes n}$ with $C(\mathcal{P}_{n})$ as in Proposition
\ref{pro:Embeddings}. We denote $i_{n(v),n}(\rho_{n(v)}^{\ast}(g_{v}))\in C(\mathcal{P}_{n})$
also by $g_{v}$ for each $v\in V$. Let $g:=\sum_{v\in V}g_{v}\in C(\mathcal{P}_{n})$.
Then $\langle g,g\rangle_{A}(y)\ge1$ for all $y\in K$. Thus $b:=\langle g,g\rangle_{A}$
is positive and invertible. Let $\xi:=gb^{-1/2}\in X^{\otimes n}$.
Then\[
\langle\xi,\xi\rangle_{A}=\langle gb^{-1/2},gb^{-1/2}\rangle_{A}=b^{-1/2}\langle g,g\rangle_{A}b^{-1/2}=1_{A}.\]
For any $y\in K$ and any $\alpha\in E^{n}$ such that $y\in K_{r(\alpha)}$,
let $x=\phi_{\alpha}(y)$. If $x\in U_{0}$, then $\Vert a\Vert-\varepsilon\le a(x)$,
and, if $x\in U_{0}^{c}$, then\[
\xi(\phi_{\alpha_{1},\dots,\alpha_{n}}(y),\dots,\phi_{\alpha_{n}}(y),y)=g(x,y)b^{-1/2}(y)=0.\]
Therefore\begin{eqnarray*}
\Vert a\Vert-\varepsilon & = & (\Vert a\Vert-\varepsilon)\langle\xi,\xi\rangle_{A}(y)\\
 & = & (\Vert a\Vert-\varepsilon)\sum_{\twolines{\alpha\in E^{n}}{y\in K_{r(\alpha)}}}|\xi(\phi_{\alpha_{1},\dots,\alpha_{n}}(y),\dots,\phi_{\alpha_{n}}(y),y)|^{2}\\
 & \le & \sum_{\twolines{\alpha\in E^{n}}{y\in K_{r(\alpha)}}}a(\phi_{\alpha}(y))|\xi(\phi_{\alpha_{1},\dots,\alpha_{n}}(y),\dots,\phi_{\alpha_{n}}(y),y)|^{2}\\
 & = & \langle\xi,a\xi\rangle_{A}(y)=S_{\xi}^{\ast}aS_{\xi}(y).\end{eqnarray*}
We also have that\[
S_{\xi}^{\ast}aS_{\xi}=\langle\xi,a\xi\rangle_{A}\le\Vert a\Vert\langle\xi,\xi\rangle_{A}=\Vert a\Vert.\]

\end{proof}
\begin{lem}
Let $\MWG$ be a Mauldin-Williams graph. Assume that $K$ is the invariant
set of the graph. For any non-zero positive element $a\in A$ and
for any $\varepsilon>0$ with $0<\varepsilon<\Vert a\Vert$, there
exists $n\in\mathbb{N}$ and $u\in X^{\otimes n}$ such that\[
\Vert u\Vert_{2}\le(\Vert a\Vert-\varepsilon)^{-1/2}\;\mbox{and}\; S_{u}^{\ast}aS_{u}=1.\]

\end{lem}
\begin{proof}
The proof is identical with the proof of \cite[Lemma 3.4]{KWp03}.
\end{proof}
\begin{lem}
Let $\MWG$ be a Mauldin-Williams graph. Suppose that the graph $G$
has no sinks and no sources and it is irreducible and not a cyclic
permutation. Assume that $K$ is the invariant set of the graph and
that $\mathcal{M}$ satisfies the open set condition in $K$. For
any $n\in\mathbb{N}$, for any $T\in\mathcal{L}(X^{\otimes n})$,
and for any $\varepsilon>0$ there exists a positive element $a\in A$
such that $a$ is $(\mathcal{M},E^{n})$-invariant,\[
\Vert\Phi(a)T\Vert^{2}\ge\Vert T\Vert^{2}-\varepsilon,\]
and $\beta^{p}(a)\beta^{q}(a)=0$ for $p,q=1,\dots,n$ with $p\ne q$.
\end{lem}
\begin{proof}
Let $n\in\mathbb{N}$, let $T\in\mathcal{L}(X^{\otimes n})$, and
let $\varepsilon>0$. Then there exists $\xi\in X^{\otimes n}$ such
that $\Vert\xi\Vert_{2}=1$ and $\Vert T\Vert^{2}\ge\Vert T\xi\Vert_{2}^{2}>\Vert T\Vert^{2}-\varepsilon$.
Hence there exists $y_{0}\in K_{v_{0}}$ for some $v_{0}\in V$ such
that\[
\Vert T\xi\Vert_{2}^{2}=\sum_{\twolines{\alpha\in E^{n}}{r(\alpha)=v_{0}}}|(T\xi)(\phi_{\alpha_{1},\dots,\alpha_{n}}(y_{0}),\dots,\phi_{\alpha_{n}}(y_{0}),y_{0})|^{2}>\Vert T\Vert^{2}-\varepsilon.\]
Then there exists an open neighborhood $U_{0}$ of $y_{0}$ in $K_{v_{0}}$
such that for any $y\in U_{0}$\[
\sum_{\twolines{\alpha\in E^{n}}{r(\alpha)=v_{0}}}|(T\xi)(\phi_{\alpha_{1},\dots,\alpha_{n}}(y),\dots,\phi_{\alpha_{n}}(y),y)|^{2}>\Vert T\Vert^{2}-\varepsilon.\]
Since $\mathcal{M}$ satisfies the open set condition in $K$, there
exists a family of non-empty sets $(V_{v})_{v\in E^{0}}$, such that
$V_{v}\subset K_{v}$, for all $v\in E^{0}$,\[
\bigcup_{\twolines{e\in E^{1}}{s(e)=v}}\phi_{e}(V_{r(e)})\subset V_{v}\;\mbox{for all }v\in E^{0},\]
\[
\mbox{and }\phi_{e}(V_{r(e)})\bigcap\phi_{f}(V_{r(f)})=\emptyset\;\textrm{\mbox{if }}e\ne f.\]
Then there exists $y_{1}\in V_{v_{0}}\cap U_{0}$ and an open neighborhood
$U_{1}$ of $y_{1}$ with $U_{1}\subset V\cap U_{0}$. Moreover, there
is some $k^{\prime}\in\mathbb{N}$ and $(e_{1},\dots,e_{k^{\prime}})\in E^{k^{\prime}}(v_{0})$
such that\[
\phi_{e_{1},\dots,e_{k^{\prime}}}(V_{r(e_{k^{\prime}})})\subset U_{1}\subset V_{v_{0}}\cap U_{0}.\]
Since the graph $G$ is not a cyclic permutation, there is a vertex
$v^{\ast}\in E^{0}$ and two edges $e^{\prime},e^{\prime\prime}\in E^{1}$
such that $e^{\prime}\ne e^{\prime\prime}$. Since the graph $G$
is irreducible, there exists a path from $r(e_{k^{\prime}})$ to $v^{\ast}$.
Hence we have a path $(e_{1},\dots,e_{k})\in E^{k}(v_{0})$ for some
$k\in\mathbb{N}$, $k\ge k^{\prime}$, such that $r(e_{k})=v^{\ast}$
and $\phi_{e_{1},\dots,e_{k}}(V_{r(e_{k})})\subset U_{1}\subset V_{v_{0}}\cap U_{0}.$
Then we can find a path $(e_{k+1},\dots,e_{k+n})\in E^{n}(v^{\ast})$
such that $e_{k+1}\ne e_{k+i}$ if $i\ne1$. To see this, let $e_{k+1}=e^{\prime}$.
If $r(e_{k+1})=v^{\ast}$, take $e_{k+2}=e^{\prime\prime}$. If $r(e_{k+1})\ne v^{\ast}$,
since $G$ has no sinks, there is an edge $e\in E^{1}$ such that
$s(e)=r(e_{k+1})$. Then $e\ne e_{k+1}$. Let $e_{k+2}=e$. If $r(e_{k+2})=v^{\ast}$,
take $e_{k+3}=e^{\prime\prime}$; if $r(e_{k+2})\ne v^{\ast}$, take
$e_{k+3}$ to be any edge such that $s(e_{k+3})=r(e_{k+2})$. Therefore
$e_{k+3}\ne e_{k+1}$. Inductively, we obtain the path $(e_{k+1},\dots,e_{k+n})\in E^{n}(v^{\ast})$
with the desired property. Then\[
\emptyset\ne\phi_{e_{1},\dots,e_{k+n}}(V_{r(e_{k+n})})\subset U_{1}\subset V_{v_{0}}\cap U_{0}.\]
There exist $y_{2}\in U_{1}$, an open neighborhood $U_{2}$ of $y_{2}$
in $K_{v_{0}}$ and a compact set $L$ such that\[
y_{2}\in U_{2}\subset L\subset\phi_{e_{1},\dots,e_{k+n}}(V_{r(e_{k+n})})\subset U_{1}\subset V_{v_{0}}\cap U_{0}.\]
Let $b\in A$ such that $0\le b\le1$, $b(y_{2})=1$ and $b|_{U_{2}^{c}}=0$.
For $\alpha\in E^{n}$ such that $r(\alpha)=v_{0}$, we have\[
\phi_{\alpha}(y_{2})\in\phi_{\alpha}(U_{2})\subseteq\phi_{\alpha}(L)\subseteq\phi_{\alpha}(V_{v_{0}}).\]
Moreover, for $\alpha,\beta\in E^{n}$ such that $r(\alpha)=r(\beta)=v_{0}$,
by the open set condition,\[
\phi_{\alpha}(L)\cap\phi_{\beta}(L)=\emptyset\;\mbox{if}\;\alpha\ne\beta.\]
We define a positive function $a$ on $K$ by the formula\[
a(x)=\begin{cases}
b(\phi_{\alpha}^{-1}(x)) & \mbox{if}\; x\in\phi_{\alpha}(L),\;\alpha\in E^{n}\;\mbox{such that}\; r(\alpha)=v_{0}\\
0 & \mbox{otherwise}.\end{cases}\]
Then $a$ is continuous on $K$, hence $a\in A$. By construction,
$a$ is $(\mathcal{M},E^{n})$-invariant.
\newcommand{\supp}{\operatorname{supp}}

Let $p\le n$ be a natural number. Let $(\alpha_{1},\dots,\alpha_{p})\in E^{p}$.
If there is no path $(\alpha_{p+1},\dots,\alpha_{n})\in E^{n-p}(r(\alpha_{p}))$
such that $r(\alpha_{n})=v_{0}$, then $\beta_{\alpha_{p}}\dots\beta_{\alpha_{1}}(a)=0$.
If there is at least one path $(\alpha_{p+1},\dots,\alpha_{n})\in E^{n-p}(r(\alpha_{p}))$
such that $r(\alpha_{n})=v_{0}$, then\[
\supp(\beta_{\alpha_{p}}\dots\beta_{\alpha_{1}}(a))\subseteq\bigcup_{\twolines{(\alpha_{p+1},\dots,\alpha_{n})\in E^{n-p}(r(\alpha_{p}))}{r(\alpha_{n})=v_{0}}}\phi_{\alpha_{p+1}\dots\alpha_{n}}(\supp(b)).\]
Since $\supp(b)\subseteq L\subset\phi_{e_{1}\dots e_{k+n}}(V_{r(e_{k+n})})$
we have that\[
\supp(\beta_{\alpha_{p}}\dots\beta_{\alpha_{1}}(a))\subseteq\bigcup_{\twolines{(\alpha_{p+1},\dots,\alpha_{n})\in E^{n-p}(r(\alpha_{p}))}{r(\alpha_{n})=v_{0}}}\phi_{\alpha_{p+1}\dots\alpha_{n}}\phi_{e_{1}\dots e_{k+n}}(V_{r(e_{k+n})}).\]
Then, for $1\le p<q\le n$ we have that\[
\supp(\beta^{p}(a))\subseteq\bigcup_{\twolines{(\alpha_{p+1},\dots,\alpha_{n})\in E^{n-p}}{r(\alpha_{n})=v_{0}}}\phi_{\alpha_{p+1}\dots\alpha_{n}}\phi_{e_{1}\dots e_{k+n}}(V_{r(e_{k+n})})\]
and\[
\supp(\beta^{q}(a))\subseteq\bigcup_{\twolines{(\alpha_{q+1},\dots,\alpha_{n})\in E^{n-q}}{r(\alpha_{n})=v_{0}}}\phi_{\alpha_{q+1}\dots\alpha_{n}}\phi_{e_{1}\dots e_{k+n}}(V_{r(e_{k+n})}).\]
Since $(n-p)+(k+1)$-th subsuffixes are different as $e_{k+1}\ne e_{k+1+(q-p)}$,
we have that $\supp(\beta^{p}(a))\cap\supp(\beta^{q}(a))=\emptyset.$
Hence $\beta^{p}(a)\beta^{q}(a)=0$.

Moreover, we have\begin{eqnarray*}
\Vert\Phi(a)T\xi\Vert_{2}^{2} & = & \sup_{y\in K}\sum_{\twolines{\alpha\in E^{n}}{y\in K_{r(\alpha)}}}\mid(a(\phi_{\alpha}(y))(T\xi)(\phi_{\alpha}(y),\dots,\phi_{\alpha_{n}}(y),y)\mid^{2}\\
 & = & \sup_{y\in L}\sum_{\twolines{\alpha\in E^{n}}{y\in K_{r(\alpha)}}}\mid(b(y))(T\xi)(\phi_{\alpha}(y),\dots,\phi_{\alpha_{n}}(y),y)\mid^{2}\\
 & \ge & \sum_{\twolines{\alpha\in E^{n}}{y_{2}\in K_{r(\alpha)}}}\mid(T\xi)(\phi_{\alpha}(y_{2}),\dots,\phi_{\alpha_{n}}(y_{2}),y_{2})\mid^{2}\\
 & > & \Vert T\Vert^{2}-\varepsilon,\end{eqnarray*}
because $y_{2}\in U_{0}$. Thus $\Vert\Phi(a)T\Vert^{2}\ge\Vert T\Vert^{2}-\varepsilon$.
\end{proof}
As in \cite{KWp03}, we let $\mathcal{F}_{n}$ be the $C^{\ast}$-subalgebra
of $\mathcal{F}_{X}$ generated by $\mathcal{K}(X^{\otimes k})$,
$k=0,1,\dots,n$ and $B_{n}$ be the $C^{\ast}$-subalgebra of $\mathcal{O}_{X}$
generated by\[
\bigcup_{i=1}^{n}\{ S_{x_{1}}\cdots S_{x_{k}}S_{y_{k}}^{\ast}\dots S_{y_{1}}^{\ast}\;:\; x_{1},\dots,x_{k},y_{1},\dots,y_{k}\in X\}\cup A.\]
We will also use the isomorphism $\varphi:\mathcal{F}_{n}\rightarrow B_{n}$
such that\[
\varphi(\theta_{x_{1}\otimes\dots\otimes x_{k},y_{1}\otimes\dots\otimes y_{k}})=S_{x_{1}}\dots S_{x_{k}}S_{y_{k}}^{\ast}\dots S_{y_{1}}^{\ast}.\]

\begin{lem}
In the above situation, let $b=c^{\ast}c$ for some $c\in\mathcal{O}_{X}^{alg}$.
We decompose $b=\sum_{j}b_{j}$ with $\gamma_{t}(b_{j})=e^{ijt}b_{j}$.
For any $\varepsilon>0$ there exists $P\in A$ with $0\le P\le I$
satisfying the following:
\begin{enumerate}
\item $Pb_{j}P=0$ $(j\ne0)$.
\item $\Vert Pb_{0}P\Vert\ge\Vert b_{0}\Vert-\varepsilon$.
\end{enumerate}
\end{lem}
\begin{proof}
The proof requires only small modifications from the proof of \cite[Lemma 3.6]{KWp03}.
\end{proof}
Having proved the equivalent of the \cite[Lemma 3.1-Lemma 3.6]{KWp03},
we obtain, using the same proof as \cite[Theorem 3.7]{KWp03}, the
corresponding result for the Mauldin-Williams graph:

\begin{thm}
Let $\MWG$ be a Mauldin-Williams graph. Suppose that the graph $G$
has no sinks and no sources, it is irreducible, and is not a cyclic
permutation. Assume that $K$ is the invariant set of the Mauldin-Williams
graph and that $\mathcal{M}$ satisfies the open set condition in
$K$. Then the associated $C^{\ast}$-algebra $\mathcal{O}_{\mathcal{M}}(K)$
is simple and purely infinite.
\end{thm}
Using the same argument as in \cite[Proposition 3.8]{KWp03} one can
show that $\mathcal{O}_{\mathcal{M}}(K)$ is separable and nuclear,
and satisfies the Universal Coefficient Theorem. Thus, by the classification
theorem of Kirchberg and Phillips \cite{Kirchberg,Phillips}, the
isomorphism class of $\mathcal{O}_{\mathcal{M}}(K)$ is completely
determined by the $K$-theory with the class of the unit.

\section{Examples}

We will compute the $K$-groups of the $C^{\ast}$-algebra associated
with a graph $G=(E^{0},E^{1},r,s)$ using the fact that $K_{1}(C^{\ast}(G))$
is isomorphic to the kernel of $1-A_{G}^{t}:\mathbb{Z}^{E^{0}}\rightarrow\mathbb{Z}^{E^{0}}$,
and $K_{0}(C^{\ast}(G))$ is isomorphic to the cokernel of the same
map, where $A_{G}$ is the vertex $E^{0}\times E^{0}$ matrix defined
by\[
A_{G}(v,w)=\#\{ e\in E^{1}\;:\; s(e)=v\;\mbox{and}\; r(e)=w\}.\]
For the $K$-groups of the Cuntz-Pimsner algebras we will use the
following six-term exact sequence due to Pimsner \cite{Pimnser} (see
also \cite[Section 4]{KWp03})

\hspace{6pt}

\begin{center}

\input{SixExactSequence.pstex_t}

\end{center}

\hspace{6pt}

First we give a condition, which is similar to the one in \cite[Section 4]{KWp03},
that implies that the associated $C^{\ast}$-algebra $\mathcal{O}_{\mathcal{M}}(K)$
is isomorphic to the $C^{\ast}$-algebra associated with the underlying
graph.

\begin{defn}
We say that a Mauldin-Williams graph $\MWG$ satisfies the \emph{graph
separation condition} in $K$ if $K$ is the invariant set of the
Mauldin-Williams graph and\[
\cograph(\phi_{e})\bigcap\cograph(\phi_{f})=\emptyset\;\mbox{if}\; e\ne f.\]

\end{defn}
\begin{prop}
Let $\MWG$ be a Mauldin-Williams graph which satisfies the graph
separation condition. Then the associated $C^{\ast}$-algebra is isomorphic
to the $C^{\ast}$-algebra associated with the underlying graph $G$.
\end{prop}
\begin{proof}
Let $E^{1}\times_{G}K=\{(e,x)\:|\: x\in K_{r(e)}\}$. Let $\mathcal{X}$
be the $C^{\ast}$-correspondence over $A=C(K)$ associated in \cite{John}
with a Mauldin-Williams graph. That is, $\mathcal{X}=C(E^{1}\times_{G}K)$
with the operations given by the formulae:\[
\xi\cdot a(e,x):=\xi(e,x)a(x),\]
\[
a\cdot\xi(e,x):=a\circ\phi_{e}(x)\xi(e,x),\]
 where $a\in C(K)$ and $\xi\in\mathcal{X}$, and\[
\langle\xi,\eta\rangle_{A}(x):=\sum_{{\genfrac{}{}{0pt}{}{e\in E^{1}}{x\in K_{r(e)}}}}\overline{\xi(e,x)}\eta(e,x),\]
 for $\xi,\eta\in\mathcal{X}$. It was proved in \cite[Theorem 2.3]{John}
that the Cuntz-Pimsner algebra associated with this $C^{\ast}$-correspondence
is isomorphic to the Cuntz-Krieger $C^{\ast}$-algebra $C^{\ast}(G)$
of the underlying graph $G$. We will show that if the Mauldin-Williams
graph satisfies the graph separation condition, then $X$ and $\mathcal{X}$
are isomorphic as $C^{\ast}$-correspondences in the sense of \cite{Muhly-Baruc-Moritaequiv}.

Let $V:X\rightarrow\mathcal{X}$ defined by the formula\[
(V\xi)(e,x)=\xi(\phi_{e}(x),x)\;\mbox{for all}\;(e,x)\in E^{1}\times_{G}K.\]
Then\begin{eqnarray*}
(V(a\cdot\xi\cdot b))(e,x) & = & (a\cdot\xi\cdot b)(\phi_{e}(x),x)\\
 & = & a(\phi_{e}(x))\xi(\phi_{e}(x),x)b(x)=(a\cdot V(\xi)\cdot b)(e,x),\end{eqnarray*}
for all $a,b\in A$ and $\xi\in X$. Also\begin{eqnarray*}
\langle V(\xi),V(\eta)\rangle_{A}(x) & = & \sum_{\twolines{e\in E^{1}}{x\in K_{r(e)}}}\overline{V(\xi)(e,x)}V(\eta)(e,x)\\
 & = & \sum_{\twolines{e\in E^{1}}{x\in K_{r(e)}}}\overline{\xi(\phi_{e}(x),x)}\eta(\phi_{e}(x),x)=\langle\xi,\eta\rangle_{A}(x),\end{eqnarray*}
for all $\xi,\eta\in X$. Finally, let $\eta\in\mathcal{X}$. For
$(y,x)\in\mathcal{G}$, since $\mathcal{M}$ satisfies the graph separation
condition, there is a unique $e\in E^{1}$ such that $y=\phi_{e}(x)$.
Define $\xi\in X$ by the formula\[
\xi(y,x)=\eta(e,x)\;\mbox{for all}\;(y,x)\in\mathcal{G}.\]
Then $\xi$ is well defined and $V(\xi)=\eta$, hence $V$ is onto.
Thus $V$ is a $C^{\ast}$-correspondence isomorphism. Therefore $\mathcal{O}_{\mathcal{M}}(K)$
is isomorphic to $C^{\ast}(G)$.
\end{proof}
\begin{example}
(The two-part dust) Let $G$ be the graph from the figure

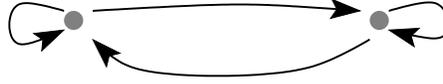
\hspace{6pt}

\begin{center}

\input{Twopartdust.pstex_t}

\end{center}

\hspace{6pt}

\noindent Let $T_{1},T_{2}\subset\mathbb{R}^{2}$ be two disjoint
compact sets such that $K_{1}$ contains the origin, $K_{2}$ contains
$(0,1)$, and $K_{1}$, $K_{2}$ are symmetric about $x$-axis. Let
$\{\phi_{i}\}_{i=1,\dots,4}$ be similarities, such that the map $\phi_{1}$
has ratio $1/2$, fixed point $(0,0)$, and rotation $30$ degrees
counterclockwise. The map $\phi_{2}$ has ratio $1/4$, fixed point
$(1,0)$, and rotation $60$ degrees clockwise. The map $\phi_{3}$
has ratio $1/2$, fixed point $(0,0)$ and rotation $90$ degrees
counterclockwise. The map $\phi_{4}$ has ratio $3/4$, fixed point
$(1,0)$, and rotation $120$ degrees clockwise (see \cite[Page 167]{Edgar}
for more details). Then the Mauldin-Williams graph $\mathcal{M}$
satisfies the graph separation condition, hence $\mathcal{O}_{\mathcal{M}}(K)$
is isomorphic to $C^{\ast}(G)$.

One can see that any Mauldin-Williams graph associated with this graph
$G$ will satisfy the graph separation condition, hence the associated
Cuntz-Pimsner algebra will be isomorphic to $C^{\ast}(G)$.
\end{example}
\;

\begin{example}
Let $G=(E^{0},E^{1},r,s)$ be the graph with the vertex matrix $A_{G}$\[
A_{G}=\left[\begin{array}{cc}
3 & 1\\
1 & 3\end{array}\right].\]
Hence $K_{1}(C^{\ast}(G))=0$ and $K_{0}(C^{\ast}(G))\simeq\mathbb{Z}/3\mathbb{Z}$.
Let $K_{1}$ and $K_{2}$ be two copies of the square of length one
with the maps $\{\phi_{i}\}_{i=1,\dots,8}$ as in the following diagram:

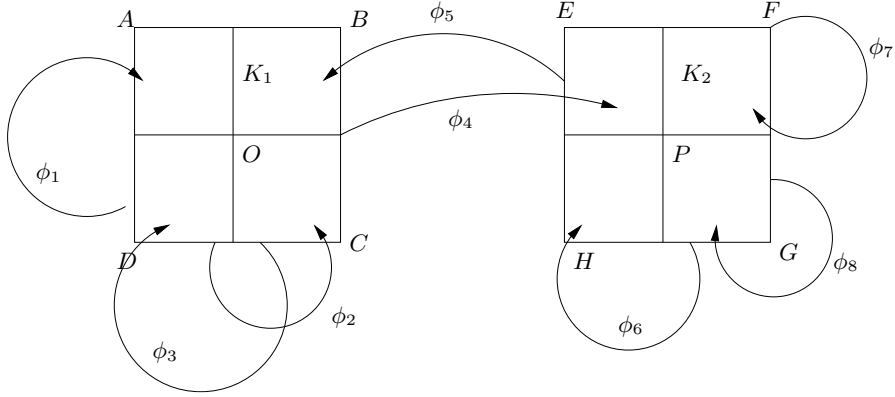
\hspace{6pt}

\begin{center}

\input{MWGraph.pstex_t}

\end{center}

\hspace{6pt}

\noindent That is $\{\phi_{i}\}_{i=1,\dots,8}$ are similarities,
$\phi_{1},\phi_{2},\phi_{3}:K_{1}\rightarrow K_{1}$, $\phi_{4}:K_{1}\rightarrow K_{2}$,
$\phi_{5}:K_{2}\rightarrow K_{1}$, $\phi_{6},\phi_{7},\phi_{8}:K_{2}\rightarrow K_{2}$,
such that $\phi_{1}(A)=A$, $\phi_{1}(C)=\phi_{2}(C)=O$, $\phi_{2}(A)=C$,
$\phi_{3}(D)=D$, $\phi_{3}(B)=O$, $\phi_{4}(A)=E$, $\phi_{4}(C)=P$,
$\phi_{5}(F)=B$, $\phi_{5}(H)=0$, $\phi_{6}(F)=P$, $\phi_{6}(H)=H$,
$\phi_{7}(F)=F$, $\phi_{7}(H)=P$, $\phi_{8}(E)=P$, $\phi_{8}(G)=G$.
Then $K=K_{1}\sqcup K_{2}$ is the invariant set of the Mauldin-Williams
graph and $B(\mathcal{M})=\{ O\}$. So $A=C(K_{1})\oplus C(K_{2})$,
$I_{X}=C_{0}(K_{1}\setminus\{ O\})\oplus C(K_{2})$. Then $K_{1}(\mathcal{O}_{\mathcal{M}}(K))\simeq0$
and $K_{0}(\mathcal{O}_{\mathcal{M}}(K))\simeq\mathbb{Z}/2\mathbb{Z}$.
Therefore $\mathrm{K}_{0}(\mathcal{O}_{\mathcal{M}}(K))$ is not isomorphic
to $\mathrm{K}_{0}(C^{\ast}(G))$. Thus $\mathcal{O}_{\mathcal{M}}(K)$
is not isomorphic to $C^{\ast}(G)$.
\end{example}

\begin{example}
(Penrose tiling) The tiles of the Penrose tiling are two triangles
such that the angles of the first triangle are equal to $\frac{\pi}{5},\frac{2\pi}{5},\frac{2\pi}{5}$
and the angles of the second triangle are equal to $\frac{3\pi}{5},\frac{\pi}{5},\frac{\pi}{5}$.
These triangles are cut from a rectangular pentagon by the diagonals
issued from a common vertex. Hence the ratio between the lengths of
the shorter and the longer sides of the triangles is equal to $\tau=\frac{1+\sqrt{5}}{2}$.
The Mauldin-Williams graph associated with the Penrose tiling is given
in the following diagram as shown in \cite{BarGrigNek}:

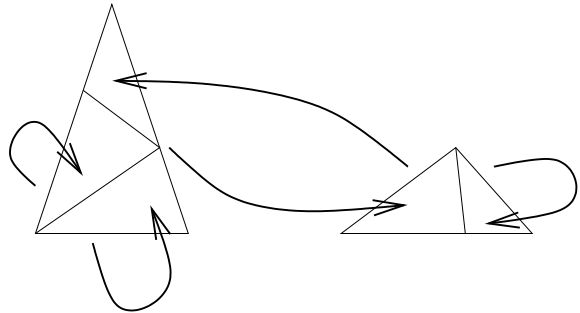
\hspace{6pt}

\begin{center}

\input{Penrose.pstex_t}

\end{center}

\hspace{6pt}

\noindent Then the invariant set $K$ of this Mauldin-Williams is
the union of the two triangles. The vertex matrix of the graph is\[
A_{G}=\left[\begin{array}{cc}
2 & 1\\
1 & 1\end{array}\right],\]
hence $K_{0}(C^{\ast}(G))=K_{1}(C^{\ast}(G))=0$. Since $K_{0}(\mathcal{O}_{\mathcal{M}}(K))\simeq\mathbb{Z}$
and $K_{1}(\mathcal{O}_{\mathcal{M}}(K))=0$, $\mathcal{O}_{\mathcal{M}}(K)$
is not isomorphic to $C^{\ast}(G)$.
\end{example}

\end{document}

%% file: SixExactSequence.pstex_t
\begin{picture}(0,0)%
\includegraphics{SixExactSequence.pstex}%
\end{picture}%
\setlength{\unitlength}{3158sp}%
\begingroup\makeatletter\ifx\SetFigFont\undefined%
\gdef\SetFigFont#1#2#3#4#5{%
  \reset@font\fontsize{#1}{#2pt}%
  \fontfamily{#3}\fontseries{#4}\fontshape{#5}%
  \selectfont}%
\fi\endgroup%
\begin{picture}(6538,1864)(1051,-1244)
\put(2551,464){\makebox(0,0)[lb]{\smash{{\SetFigFont{10}{12.0}{\familydefault}{\mddefault}{\updefault}{\color[rgb]{0,0,0}$id-[X]$}%
}}}}
\put(5176,-1186){\makebox(0,0)[lb]{\smash{{\SetFigFont{10}{12.0}{\familydefault}{\mddefault}{\updefault}{\color[rgb]{0,0,0}$id-[X]$}%
}}}}
\put(5026,464){\makebox(0,0)[lb]{\smash{{\SetFigFont{10}{12.0}{\familydefault}{\mddefault}{\updefault}{\color[rgb]{0,0,0}$i_\ast$}%
}}}}
\put(2776,-1186){\makebox(0,0)[lb]{\smash{{\SetFigFont{10}{12.0}{\familydefault}{\mddefault}{\updefault}{\color[rgb]{0,0,0}$i_\ast$}%
}}}}
\put(6976,-286){\makebox(0,0)[lb]{\smash{{\SetFigFont{10}{12.0}{\familydefault}{\mddefault}{\updefault}{\color[rgb]{0,0,0}$\delta_0$}%
}}}}
\put(1201,-361){\makebox(0,0)[lb]{\smash{{\SetFigFont{10}{12.0}{\familydefault}{\mddefault}{\updefault}{\color[rgb]{0,0,0}$\delta_1$}%
}}}}
\put(1201,239){\makebox(0,0)[lb]{\smash{{\SetFigFont{10}{12.0}{\familydefault}{\mddefault}{\updefault}{\color[rgb]{0,0,0}$K_0(I_X)$}%
}}}}
\put(6451,-961){\makebox(0,0)[lb]{\smash{{\SetFigFont{10}{12.0}{\familydefault}{\mddefault}{\updefault}{\color[rgb]{0,0,0}${K}_1(I_X)$}%
}}}}
\put(1051,-961){\makebox(0,0)[lb]{\smash{{\SetFigFont{10}{12.0}{\familydefault}{\mddefault}{\updefault}{\color[rgb]{0,0,0}$K_1(\mathcal{O}_\mathcal{M}(K))$}%
}}}}
\put(3901,-961){\makebox(0,0)[lb]{\smash{{\SetFigFont{10}{12.0}{\familydefault}{\mddefault}{\updefault}{\color[rgb]{0,0,0}$K_1(A)$}%
}}}}
\put(3901,239){\makebox(0,0)[lb]{\smash{{\SetFigFont{10}{12.0}{\familydefault}{\mddefault}{\updefault}{\color[rgb]{0,0,0}$K_0(A)$}%
}}}}
\put(6151,239){\makebox(0,0)[lb]{\smash{{\SetFigFont{10}{12.0}{\familydefault}{\mddefault}{\updefault}{\color[rgb]{0,0,0}$K_0(\mathcal{O}_{\mathcal{M}}(K))$}%
}}}}
\end{picture}%

%% file: Twopartdust.pstex_t
\begin{picture}(0,0)%
\includegraphics{Twopartdust.pstex}%
\end{picture}%
\setlength{\unitlength}{3158sp}%
\begingroup\makeatletter\ifx\SetFigFont\undefined%
\gdef\SetFigFont#1#2#3#4#5{%
  \reset@font\fontsize{#1}{#2pt}%
  \fontfamily{#3}\fontseries{#4}\fontshape{#5}%
  \selectfont}%
\fi\endgroup%
\begin{picture}(3490,646)(675,-219)
\end{picture}%

%% file: MWGraph.pstex_t
\begin{picture}(0,0)%
\includegraphics{MWGraph.pstex}%
\end{picture}%
\setlength{\unitlength}{2960sp}%
\begingroup\makeatletter\ifx\SetFigFont\undefined%
\gdef\SetFigFont#1#2#3#4#5{%
  \reset@font\fontsize{#1}{#2pt}%
  \fontfamily{#3}\fontseries{#4}\fontshape{#5}%
  \selectfont}%
\fi\endgroup%
\begin{picture}(7223,3292)(128,-2822)
\put(1051,239){\makebox(0,0)[lb]{\smash{\SetFigFont{9}{10.8}{\familydefault}{\mddefault}{\updefault}{\color[rgb]{0,0,0}$A$}%
}}}
\put(3001,239){\makebox(0,0)[lb]{\smash{\SetFigFont{9}{10.8}{\familydefault}{\mddefault}{\updefault}{\color[rgb]{0,0,0}$B$}%
}}}
\put(3001,-1636){\makebox(0,0)[lb]{\smash{\SetFigFont{9}{10.8}{\familydefault}{\mddefault}{\updefault}{\color[rgb]{0,0,0}$C$}%
}}}
\put(1051,-1786){\makebox(0,0)[lb]{\smash{\SetFigFont{9}{10.8}{\familydefault}{\mddefault}{\updefault}{\color[rgb]{0,0,0}$D$}%
}}}
\put(2101,-886){\makebox(0,0)[lb]{\smash{\SetFigFont{9}{10.8}{\familydefault}{\mddefault}{\updefault}{\color[rgb]{0,0,0}$O$}%
}}}
\put(2101,-211){\makebox(0,0)[lb]{\smash{\SetFigFont{9}{10.8}{\familydefault}{\mddefault}{\updefault}{\color[rgb]{0,0,0}$K_1$}%
}}}
\put(5776,-211){\makebox(0,0)[lb]{\smash{\SetFigFont{9}{10.8}{\familydefault}{\mddefault}{\updefault}{\color[rgb]{0,0,0}$K_2$}%
}}}
\put(376,-1036){\makebox(0,0)[lb]{\smash{\SetFigFont{9}{10.8}{\familydefault}{\mddefault}{\updefault}{\color[rgb]{0,0,0}$\phi_1$}%
}}}
\put(2851,-2236){\makebox(0,0)[lb]{\smash{\SetFigFont{9}{10.8}{\familydefault}{\mddefault}{\updefault}{\color[rgb]{0,0,0}$\phi_2$}%
}}}
\put(1351,-2536){\makebox(0,0)[lb]{\smash{\SetFigFont{9}{10.8}{\familydefault}{\mddefault}{\updefault}{\color[rgb]{0,0,0}$\phi_3$}%
}}}
\put(4726,314){\makebox(0,0)[lb]{\smash{\SetFigFont{9}{10.8}{\familydefault}{\mddefault}{\updefault}{\color[rgb]{0,0,0}$E$}%
}}}
\put(6451,314){\makebox(0,0)[lb]{\smash{\SetFigFont{9}{10.8}{\familydefault}{\mddefault}{\updefault}{\color[rgb]{0,0,0}$F$}%
}}}
\put(6601,-1711){\makebox(0,0)[lb]{\smash{\SetFigFont{9}{10.8}{\familydefault}{\mddefault}{\updefault}{\color[rgb]{0,0,0}$G$}%
}}}
\put(4876,-1786){\makebox(0,0)[lb]{\smash{\SetFigFont{9}{10.8}{\familydefault}{\mddefault}{\updefault}{\color[rgb]{0,0,0}$H$}%
}}}
\put(3826,-586){\makebox(0,0)[lb]{\smash{\SetFigFont{9}{10.8}{\familydefault}{\mddefault}{\updefault}{\color[rgb]{0,0,0}$\phi_4$}%
}}}
\put(3676,314){\makebox(0,0)[lb]{\smash{\SetFigFont{9}{10.8}{\familydefault}{\mddefault}{\updefault}{\color[rgb]{0,0,0}$\phi_5$}%
}}}
\put(5251,-2311){\makebox(0,0)[lb]{\smash{\SetFigFont{9}{10.8}{\familydefault}{\mddefault}{\updefault}{\color[rgb]{0,0,0}$\phi_6$}%
}}}
\put(7351, 14){\makebox(0,0)[lb]{\smash{\SetFigFont{9}{10.8}{\familydefault}{\mddefault}{\updefault}{\color[rgb]{0,0,0}$\phi_7$}%
}}}
\put(7051,-1786){\makebox(0,0)[lb]{\smash{\SetFigFont{9}{10.8}{\familydefault}{\mddefault}{\updefault}{\color[rgb]{0,0,0}$\phi_8$}%
}}}
\put(5701,-886){\makebox(0,0)[lb]{\smash{\SetFigFont{9}{10.8}{\familydefault}{\mddefault}{\updefault}{\color[rgb]{0,0,0}$P$}%
}}}
\end{picture}

%% file: Penrose.pstex_t
\begin{picture}(0,0)%
\includegraphics{Penrose.pstex}%
\end{picture}%
\setlength{\unitlength}{3158sp}%
\begingroup\makeatletter\ifx\SetFigFont\undefined%
\gdef\SetFigFont#1#2#3#4#5{%
  \reset@font\fontsize{#1}{#2pt}%
  \fontfamily{#3}\fontseries{#4}\fontshape{#5}%
  \selectfont}%
\fi\endgroup%
\begin{picture}(4491,2439)(979,-2188)
\end{picture}%